\documentclass[11pt]{gtart}
\usepackage{amsmath, amssymb, graphicx, epsfig, graphs, verbatim, psfrag}
\usepackage{epic}
\newtheorem{thm}{Theorem}[section]

\newtheorem{cor}[thm]{Corollary}
\newtheorem{lem}[thm]{Lemma}
\newtheorem{prop}[thm]{Proposition}

\theoremstyle{definition}
\newtheorem{defn}[thm]{Definition}

\newtheorem{rem}[thm]{Remark}

\numberwithin{equation}{section}

\newcommand{\bfz}{{\mathbb {Z}}}
\newcommand{\bfq}{{\mathbb {Q}}}

\newcommand{\oa}{{\overline {{{a}}}}}
\newcommand{\ox}{{\overline {{{x}}}}}
\newcommand{\oy}{{\overline {{{y}}}}}

\newcommand{\s}{\mathbf s}

\renewcommand{\u}{\mathbf u}
\renewcommand{\t}{\mathbf t}
\newcommand{\z}{\mathbf z}

\newcommand{\hf}{{{\widehat {HF}}}}

\newcommand{\mubar}{{\overline {\mu}}}

\begin{document}

\title{On the $\mubar$ invariant of rational surface singularities}

\author{Andr\'{a}s I. Stipsicz}
\address{R\'enyi Institute of Mathematics\\
H-1053 Budapest, Re\'altanoda utca 13--15, Hungary and \\
Department of Mathematics, Columbia University, NY, New York, 10027}
\email{stipsicz@math-inst.hu, stipsicz@math.columbia.edu}

\begin{abstract}
  We show that for rational surface singularities with odd determinant
  the $\mubar$ invariant defined by W. Neumann is an obstruction for
  the link of the singularity to  bound a rational homology 4--ball. We
  identify the $\mubar$ invariant with the corresponding correction
  term in Heegaard Floer theory.
\end{abstract}
\maketitle

\section{Introduction}
\label{sec:one}

Smoothings of surface singularities play a prominent role in
constructing new and interesting smooth (and symplectic) 4--manifolds.
It is of particular interest when the singularity at hand admits a
rational homology 4--ball smoothing. Such smoothings led to the
discovery of the rational blow--down procedure \cite{FSrat, Prat},
which in turn provided a major tool for finding exotic 4--manifolds.
Restrictions for a singularity to admit rational homology 4--ball
smoothing have been found recently in \cite{SSW}.

A topological obstruction for a $\bfz _2$--homology 3--sphere (that
is, a 3--manifold $Y$ with $H_*(Y; \bfz _2)=H_* (S^3; \bfz _2)$) to
bound a spin rational homology 4--ball is its $\mu$--invariant,
defined modulo 16. An integral lift $\mubar$ of $\mu$ has been defined
by Neumann in \cite{neumann1} (cf. also \cite{savel}) for plumbed
$\bfz _2$--homology 3--spheres, but it was unclear whether this
integer valued invariant obstructs the 3--manifold to bound a spin
rational homology 4--ball. Special cases, like Seifert fibered
3--manifolds, have been considered by Saveliev \cite{savel}.  More
recently, based on work of Ozsv\'ath and Szab\'o \cite{OSzabs, OSzF1,
  OSzF2} the \emph{correction term} of spin$^c$ 3--manifolds (stemming
from gradings on the Ozsv\'ath--Szab\'o homology groups) provided
further obstructions.  For applications of these invariants along
similar lines see \cite{JabukaJosh}.

In fact, in \cite{neumann1} the $\mubar$--invariant is defined for any
spin rational homology 3--sphere which can be given by plumbing
spheres along a tree (i.e., the assumption on the parity of the
determinant of the plumbing graph can be relaxed). By identifying
$\mubar$ of a spin 3--manifold $(Y, \s)$ which is a link of a rational
surface singularity with the appropriate correction term, we show

\begin{thm} \label{t:obst}
Suppose that $Y_{\Gamma}$ is given as a plumbing of spheres along a
negative definite tree $\Gamma$, defining a rational surface
singularity.
\begin{itemize}
\item For a spin structure $\s \in Spin(Y_{\Gamma })$ the invariant
$\mubar (Y_{\Gamma} , \s )\in \bfz $ is an obstruction for the
existence of a spin$^c$ rational homology 4--ball $(X, \t )$ with boundary
$(Y_{\Gamma }, \s )$.

\item If $\Pi _{\s \in Spin (Y_{\Gamma})} \mubar ( Y_{\Gamma }, \s
)\neq 0$ then the rational singularity does not bound a spin rational
homology 4--ball.

\item Specifically, if $\det \Gamma$ is odd and $\mubar (Y_{\Gamma
})\neq 0$ then $Y_{\Gamma }$ is not the boundary of a rational
homology 4--ball.  Consequently the corresponding singularity does not
admit rational homology 4--ball smoothing.
\end{itemize}
\end{thm}

\begin{cor}
Suppose that $S_{\Gamma}$ is a normal surface singularity with $\det
\Gamma$ odd. If $\mubar (Y_{\Gamma })\neq 0$ then $S_{\Gamma }$ 
does not admit a rational homology 4--ball smoothing.
\end{cor}
\begin{proof}
  If $S_{\Gamma }$ is not a rational singularity then it does not
  admit rational homology 4--ball smoothing. If $\det \Gamma$ is odd,
  then for rational surface singularities Theorem~\ref{t:obst}
  concludes the proof.
\end{proof}
We hope that this obstruction will be useful in completing the
characteriztion of surface singularities with rational homology
4--ball smoothing, along the line initiated in \cite{SSW}. 

\begin{rem}
{\rm The assumption on the parity of $\det \Gamma$ cannot be relaxed
in general, since for example the singularity with resolution graph
having a single vertex of weight $(-4)$ has two spin structures with
$\mubar$ invariants $-3$ and $+1$, but the link of the singularity is
the boundary of a rational homology 4--ball: the complement of a
quadric in the complex projective plane. In fact, this rational
homology 4--ball can be given as a smoothing of the singularity.
In accordance with Theorem~\ref{t:obst} the spin structures on the link
of the singularity do not extend to the rational homology 4--ball.}
\end{rem}

As it was indicated earlier, the proof of Theorem~\ref{t:obst} above
rests on the following, more technical statement. Here the invariant
$d(Y, \s)$ of a spin$^c$ 3--manifold $(Y, \s)$ is the correction term
in Heegaard Floer theory.  (For more on Heegaard Floer theory see
Section~\ref{sec:four}.)

\begin{thm}\label{t:main}
Suppose that $\Gamma $ is a negative definite plumbing tree of
spheres, giving rise to a rational surface singularity. Let $\s $ be a
given spin structure on the associated 3--manifold $Y_{\Gamma}$.  Then
$\mubar ( Y_{\Gamma}, \s) =-4 d(Y_{\Gamma }, \s)$.
\end{thm}

\section{The $\mu $ and $\mubar$ invariants}
\label{sec:two}

Suppose that $Y$ is a rational homology 3--sphere, and the rank $\vert
H_1 \vert $ of its first homology is odd. Then $H_1(Y; \bfz _2)=H^1(Y;
\bfz _2)=0$, hence $Y$ admits a unique spin structure. Consider a spin
4--manifold $X$ with $\partial X=Y$.  The classical definition of
Rokhlin's $\mu$--invariant is 
\[
\mu (Y) \equiv \sigma (X)\ \ \mod 16 ,
\]
where $\sigma (X)$ is the signature of the 4--manifold $X$. The
invariance of this quantity is a simple consequence of Rokhlin's
famous result on the divisibility of the signature of a closed spin
4--manifold by 16. (If $Y$ is an integral homology sphere, that is,
$H_1(Y; \bfz )=0$ also holds, then the signature $\sigma (X)$ of a
spin 4--manifold $X$ with $\partial X=Y$ is divisible by 8, and in
this case sometimes Rokhlin's invariant is defined as $\frac{\sigma
(X)}{8} \in \bfz _2$.)

It is not hard to see that if $X$ is a spin rational homology 4--ball
(i.e., $H_*(X; \bfq )=H_* (D^4; \bfq )$) with $\partial X =Y$ and
$H_1(Y; \bfz _2)=0$ then $\mu (Y)=0$. Consequently, the
$\mu$--invariant of a $\bfz_2$--homology sphere $Y$ provides an
obstruction for $Y$ to bound a spin rational homology 4--ball.  (The
spin assumption on $X$ is important, since for example the Brieskorn
sphere $\Sigma (2,3,7)$ has $\mu =1$ and bounds a nonspin rational
homology 4--ball, cf. \cite{FSmu}.) Since $\mu$ is defined only mod
16, it is typically less effective than an integer valued invariant.
Interest in integral lifts (or related obstructions) was motivated
also by a result of Galewski and Stern \cite{GaS} about higher
dimensional (simplicial) triangulation theory.

In \cite{neumann1} Walter Neumann defined a lift $\mubar \in \bfz $ of
$\mu $ for spin 3--manifolds given by the plumbing construction along
a weighted tree. Before giving the definition of this invariant we
shortly review a few basic facts about plumbing trees. For a general
reference see \cite{neumann1}.

Suppose that $\Gamma $ is a weighted tree with nonzero determinant.
Let $X_{\Gamma }$ denote the 4--manifold defined by plumbing disk
bundles over spheres according to the weighted tree $\Gamma$, and
define $Y_{\Gamma }$ as $\partial X_{\Gamma}$. As it is described in
\cite{neumann1}, the mod 2 homology $H_1(Y_{\Gamma} ; \bfz _2)$ can be
determined by a simple algorithm, which we outline below. Consider a
leaf $v$ of $\Gamma$, connected to the vertex $w$.

\begin{itemize}
\item {\bf Move 1}:  If the weight on $v$ is even, then erase $v$ and $w$
from $\Gamma$.
\item {\bf Move 2}:  If the weight of $v$ is odd, then erase $v$ and
change the parity of the weight on $w$.
\end{itemize}

This procedure stops once we reach a graph $\Gamma '$ with no edges.
Suppose that $\Gamma '$ contains $p$ vertices, $q$ of them with even
weights.

\begin{lem}
The dimension of the vector space $H_1(Y_{\Gamma }; \bfz _2)$
over $\bfz _2$ is equal to $q$.
\end{lem}
\begin{proof}
Denote the set of vertices of the given weighted plumbing tree
$\Gamma$ with nonzero determinant by $V=\{ v_1, \ldots , v_n\}$. It is
known (cf. \cite[Proposition~5.3.11]{GS}) that the homology group
$H_1(Y _{\Gamma}; \bfz )$ admits a presentation by taking elements of
$V$ as generators, and equations
\[
n_i\cdot v_i + \sum _{j\neq i} \langle v_j , v_i \rangle \cdot v_j=0
\]
as relations ($i=1,\ldots , n$), with the convention that $n_i$ is the
weight on $v_i$, and $\langle v_j, v_i \rangle $ is one or zero
depending on whether $v_j$ and $v_i$ (as vertices of the tree
$\Gamma$) are connected or not. These relations follow easily from the
existence of Seifert surfaces for the components of the surgery link.
The mod 2 reduction of the relations (with the same generators)
provide a presentation for $H_1(Y_{\Gamma }; \bfz _2)$. Now the moves
for simplifying the graph (until it becomes a disjoint union of some
vertices) obviously correspond to base changes and expressions of
generators in terms of others.  Indeed, when {\bf Move 1} applies to
$v$ and $w$ then the relation for $v$ shows $w=0$, while the relation
for $w$ expresses $v$ in terms of the other neighbours of $w$. In the
situation of {\bf Move 2} the relation for $v$ simply asserts that
$v=w$ (mod 2). From this observation the statement easily follows: a
single point with odd weight gives rise to a 3--manifold with
vanishing first mod 2 homology, while with even weight the first mod 2
homology is 1-dimensional.
\end{proof} 

Recall that an oriented 3--manifold $Y$ always admits a spin
structure, and the space of spin structures is parametrized by the
first mod 2 cohomology $H^1 (Y; \bfz _2)(\cong H_1(Y; \bfz _2))$ of
$Y$. A convenient parametrization of the set of spin structures on the
rational homology 3--sphere $Y_{\Gamma}$ is given as follows. First we
define a set of subsets of the vertex set for every plumbing graph
$\Gamma$. We start with a graph $\Gamma '$ having no edges: in that
case consider the subsets of the vertices which contain all vertices
with odd weights.  Every such subset will give rise to a unique subset
$S\subset V$ for the original graph $\Gamma$ as follows.  We describe
the change of $S$ under one step in the process giving $\Gamma '$ from
$\Gamma$. Suppose that $\Gamma ' $ is given by {\bf Move 1} from
$\Gamma$ (via erasing $v=v_i$ and $w=v_j$), and a set $S'\subset V'$
is specified for $\Gamma '$.  Now we define the set $S\subset V$ by
taking it to be equal to $S'$ or $S'\cup \{ v_i \}$ according as the
number of indices in $S'$ adjacent to $w=v_j$ have the same parity as
$n_j$ or $n_j-1$.  If $\Gamma '$ is derived from $\Gamma$ by {\bf Move
  2} (via erasing $v_i$) then let $S$ be equal to $S'$ or $S'\cup \{
v_i \}$ depending on whether $v_j$ was in $S'$ or not.  It is not hard
to see from this algorithm that if $v_i,v_j\in S$ then $v_i$ and $
v_j$ are not connected by an egde in $\Gamma$.

Suppose now that $S\subset V$ is a subset defined as above.  Consider
the submanifold $\Sigma _S\subset X_{\Gamma}$ defined as the union of
the spheres corresponding to the vertices in $S$. Notice that since by
construction $S$ does not contain adjacent vertices, the above surface
is a disjoint union of embedded spheres.  Let $c _{S}\in H^2(X_{\Gamma
}; \bfz )$ denote the Poincar\'e dual of $\Sigma _{S}$.  The inductive
definition (and the starting condition) shows that $c_{S}$ is a
\emph{characteristic element}, that is, for every surface $\Sigma
_v\subset X_{\Gamma}$ defined by a vertex $v$ we have
\[
c_S (\Sigma _v) \equiv n_v \ \ \mod 2 .
\]
On the simply connected 4--manifold $X_{\Gamma }$ a characteristic
cohomology class uniquely specifies a spin$^c$ structure $\t _S$,
which restricts to a spin$^c$ structure $\s _S$ on the boundary
$Y_{\Gamma}$. Since $PD(c_S)=\bigcup _v \Sigma _v = \Sigma _S$ is in
$H_2(X_{\Gamma}; \bfz )$, on the boundary the spin$^c$ structure $\s
_S=\t _S\vert _{\partial X_{\Gamma}}$ has vanishing first Chern class,
therefore it is a spin structure on $Y_{\Gamma}$. Hence every subset
$S$ constructed above defines a spin structure $\s _S$ on
$Y_{\Gamma}$; the set $S$ is called the \emph{Wu set} of the
corresponding spin structure. Since this construction provides a spin
structure on the complement $X-\Sigma _S$, it is obvious that two
different sets induce different spin structures: if $S_1$ and $S_2$
differ on the vertex $v$ of even weight (in the disconnected graph our
construction started with) then only the spin structure corresponding
to the Wu set not containing $v$ will extend to the cobordism we get
by the appropriate handle attachment along $v$. In conclusion, we get
an identification of $H_1(Y_{\Gamma }; \bfz )(\cong H^1 (Y_{\Gamma };
\bfz ))$ with the set of spin structures on $Y_{\Gamma}$: take the
characteristic function of $S$ on the starting disconnected graph
$\Gamma '$ (which by the above said determines $S$), and associate to
it the corresponding first mod 2 cohomology class. Now the definition
of the $\mubar$ invariant of Neumann (cf. also \cite{neumann1}) is as
follows.

\begin{defn}
  For a spin structure $\s$ on $Y_{\Gamma }$ consider the
  corresponding Wu set $S$ and embedded Wu surface $\Sigma _S\subset
  X_{\Gamma}$. Define $\mubar (Y_{\Gamma }, \s)\in \bfz $ as the
  difference
\[
\mubar (Y_{\Gamma }, \s )=\sigma (X_{\Gamma})-[\Sigma _S]^2 . 
\]
\end{defn}
By applying the handle calculus developed in \cite{neumann2} together
with the Wu set $S$, the proof of the following statement easily
follows.
\begin{prop}$($\cite[Theorem~4.1]{neumann1}$)$
The quantity $\mubar (Y_{\Gamma }, \s) $ is an invariant of the spin
3--manifold $(Y_{\Gamma }, \s )$ and is independent of the choices
made in the definition. \qed
\end{prop}

\section{Rational singularities}
\label{sec:three}
Consider the plumbing tree $\Gamma$ and suppose that $\Gamma $ is
negative definite. According to a classical result of Grauert
\cite{grau}, for any negative definite plumbing graph there exists a
normal surface singularity such that the plumbing along the given
graph is diffeomorphic to a resolution of the singularity.

\begin{defn}
  A normal surface singularity $S_{\Gamma}$ is called \emph{rational}
  if its geometric genus $p_g=0$.  A negative definite plumbing graph
  $\Gamma $ is \emph{rational} if there is a rational singularity
  $S_{\Gamma}$ with resolution diffeomorphic to $X_{\Gamma}$.
\end{defn}

Although the singularity corresponding to a plumbing graph might not
be unique, it is known that rationality is a topological property and
can be fairly easily read off from the plumbing graph through Laufer's
algorithm.  Namely, consider the homology class
\[
K_0=\sum _{v\in \Gamma } [\Sigma _v] \in H_2 (X_{\Gamma}; \bfz ).
\]
In the $i^{th}$ step, consider the product $K_i \cdot \Sigma _{v_j}=
\langle PD(K_i), [\Sigma _{v_j}]\rangle$. If it is at least 2 then the
algorithm stops and the singularity is not rational. If the product is
nonpositive, move to the next vertex.  Finally, if the product is 1
for some $v\in \Gamma$, then replace $K_i$ with $K_{i+1}=K_i +[\Sigma
_v]$ and start checking the value of the product for all vertices of
$\Gamma$ again. If all products are nonpositive, the algorithm stops
and the graph gives rise to a rational singularity.

\begin{lem} \label{l:rat}
A rational plumbing graph is always a (negative definite) tree of
spheres, and the link is a rational homology 3--sphere. In addition,
for any vertex $v_i \in \Gamma$ the sum of its weight $n_i$ and the
number $d_i$ of its neighbours is at most 1. \qed
\end{lem}

Notice that in a rational graph a vertex with weight $(-1)$ has degree
$d\leq 2$, hence can be blown down by keeping $\Gamma$ a plumbing
tree.  For this reason, we might assume that $n_i \leq -2 $ for all
vertices $v_i\in \Gamma$.

\section{Heegaard Floer groups}
\label{sec:four}

In \cite{OSzF1, OSzF2} a set of very powerful invariants, the
Ozsv\'ath--Szab\'o homology groups $\hf (Y, \s), HF^{\pm }(Y, \s)$ and
$HF^{\infty }(Y, \s)$ of a spin$^c$ 3--manifold $(Y, \s )$ were
introduced. In the following we will use these groups and relations
among them; for a more thorough introduction see \cite{OSzF1, OSzF2,
OSzII}. Recall that a rational homology 3--sphere $Y$ is an
\emph{$L$--space} if $\hf (Y, \s )=\bfz _2$ for every spin$^c$
structure $\s \in Spin^c (Y)$.  (In the version of the theory we are
about to apply, we use $\bfz _2$--coefficients.)  In this case we can
label the unique nonzero element of $\hf (Y, \s)$ by the corresponding
spin$^c$ structure $\s$.  Recall also that for a rational homology
3--sphere $Y$ the groups are equipped with a natural
$\bfq$--grading. The grading of the unique nontrivial element of $\hf
(Y, \s )$ for an $L$--space $Y$ is called the \emph{correction term}
$d(Y, \s )$ of the spin$^c$ 3--manifold $(Y, \s )$. For the proof of
the next proposition, see for example \cite[Theorem~2.3]{JN}.

\begin{prop}\label{p:van}
Suppose that $d(Y, \s)\neq 0$. Then there is no spin$^c$ rational
homology 4--ball $(X, \t)$ with $\partial (X, \t)=(Y, \s )$. \qed
\end{prop} 

\begin{prop}\label{p:ptlan}
Suppose that $\det \Gamma $ is odd. If $d(Y_{\Gamma }, \s )\neq 0$ for 
the unique spin structure $\s$ then $Y_{\Gamma }$ does not bound any 
rational homology 4--ball.
\end{prop}
\begin{proof}
Suppose that $Y_{\Gamma}=\partial X$ for a rational homology 4--ball
$X$. Let $\varphi \colon Y_{\Gamma }\to X$ denote the embedding of the
boundary, inducing the map $\varphi _*$ on homology.  Since $|
H_1(Y_{\Gamma }; \bfz )|$ is odd, the size of the subgroup Im$\
\varphi _*$ is also odd. This implies that an odd number of spin$^c$
structures in $Spin^c(Y_{\Gamma })$ extend to $X$.  Since $\s \in Spin^c
(Y_{\Gamma })$ and its conjugate ${\overline {\s }}$ extend at the
same time, we conclude that the spin structure $\s={\overline {\s }}$
of $Y_{\Gamma }$ extends to $X$ as a spin$^c$ structure, therefore
Proposition~\ref{p:van} concludes the proof.
\end{proof}

A relation between the singularity's holomorphic structure and its
Heegaard Floer theoretic behaviour was found by A. N\'emethi:
\begin{thm}[N\'emethi, \cite{nemethi}]
Suppose that the negative definite plumbing tree $\Gamma$ gives rise
to a rational singularity. Then $Y_{\Gamma}$ is an $L$--space. \qed
\end{thm}

\section{A relation between $\mubar(Y_{\Gamma },\s )$ and $d(Y _{\Gamma},\s )$}
\label{sec:five}

The proof of our main result about the $\mubar$--invariant relies on
the identification of it with the appropriate multiple of the
$d$--invariant of the spin 3--manifold at hand.

\begin{proof}[Proof of Theorem~\ref{t:main}]
  Let $\Gamma$ be a given negative definite rational plumbing tree
  with a spin structure $\s$ (represented by its Wu set $S\subset V$).
  Let $m_{\Gamma , S}$ denote the number of those vertices $v_i\in
  \Gamma $ which are not in $S$ but $-n_i$ of the neighbours of $v_i$
  are in $S$.  (Notice that by the rationality of $\Gamma$ this means
  that $v_i$ has $-n_i$ or $-n_i+1$ neighbours and either all or all
  but one neighbours are in $S$.)

The proof of the theorem will proceed by induction on $m_{\Gamma , S}$.
Let us start with the easy case when $m_{\Gamma ,S}=0$, that is, for any
vertex $v_i$ in $\Gamma$ we have 
\begin{equation}\label{e:kicsi}
c_S(\Sigma _{v_i})< -n_i.
\end{equation}
For $v_i \in S$ we have $c_S(\Sigma _{v_i})=n_i$, while if $v_i $ is
not in $S$ then $c_S(\Sigma _{v_i})$ is the number of neighbours of
$v_i$ which are in $S$. In particular, $0\leq c_S(\Sigma _{v_i})\leq
d_i$ holds for all $v_i$ not in $S$.  Since $c_S$ is characteristic,
Inequality~\eqref{e:kicsi} actually means that $c_S(\Sigma _{v_i})\leq
-n _i-2$.  In conclusion, $c_S$ satisfies $n_i\leq c_S(\Sigma
_{v_i})\leq -n_i -2$ for all vertices, hence $c_S$ is a \emph{terminal
vector} in the sense of \cite{OSzplum}. By subtracting twice the
Poincar\'e duals of the homology classes represented by surfaces
corresponding to vertices in $S$, eventually we get a path back to a
vector $K\in H^2(X_{\Gamma}; \bfz )$ which satisfies $K(\Sigma
_{v_i})= -n_i$ for $v_i\in S$ and $K(\Sigma _{v_i})\geq -d_i\geq
n_i+2$ if $v_i$ is not in $S$. This means that $K$ is an \emph{initial
vector}, hence $c_S$ is in a \emph{full path} (again, in the
terminology of \cite{OSzplum}). By the identification of
\cite{nemethi} this implies that $c_S$ gives rise to a Heegaard Floer
homology element in $\hf (Y, \s)$ of degree $\frac{1}{4}(c^2_S-3\sigma
(X_{\Gamma})-2\chi (X_{\Gamma}))$.  (Here, as costumary in Heegaard
Floer theory, $\chi (X_{\Gamma })$ is understood as the Euler
characteristic of the cobordism we get from $S^3$ to $Y_{\Gamma }$ by
deleting a point from $X_{\Gamma}$.)  Since $Y_{\Gamma}$ is an
$L$--space, this degree must be equal to $d(Y, \s)$.  On the other
hand, since $\Gamma$ is negative definite, $\chi (X_{\Gamma})= -\sigma
(X_{\Gamma})$, hence the above formula for the degree shows that
$-\mubar(Y_{\Gamma }, \s)=c^2_S-\sigma (X_{\Gamma})$ is equal to
$4d(Y_{\Gamma}, \s )$.

Next we assume that the statement is proved for graphs $(\Gamma , S )$
with $m_{\Gamma , S}\leq m-1$. In the inductive step we will utilize
the exact triangle for Heegaard Floer homologies, proved for a surgery
triple, see \cite{OSzF2, OSzI}.  To this end, fix a graph $\Gamma$
with Wu set $S$ and corresponding spin structure $\s \in Spin
(Y_{\Gamma })$ having $m_{\Gamma , S}=m>0$ and let $v$ denote a vertex
with $-n_i$ neighbours in $S$.  (Consequently $v$ is not in $S$.)
Consider the following plumbing graphs (with spin structures specified
by the various Wu sets):
\begin{itemize}
\item Let $\Gamma ', \Gamma ''$ denote the same graphs as $\Gamma$ with
 the alteration of the framing on $v$ from $n_i$ to $n_i-2$ and
 $n_i-4$, resp.  It is easy to see that $S$ still provides Wu sets
 $S', S''$ (and hence spin structures $\s ', \s ''$) for $\Gamma '$
 and $\Gamma ''$.  Notice that $m_{\Gamma ', S'}=m_{\Gamma '',
 S''}=m_{\Gamma, S} -1$.  In addition, since $v$ was not in the Wu set
 $S$, we see directly that $\mubar (Y_{\Gamma }, \s)=\mubar (Y_{\Gamma
 '}, \s ')=\mubar (Y_{\Gamma ''},\s '')$. Laufer's algorithm shows
 that $\Gamma ', \Gamma ''$   are also rational.
 \item Let $\Gamma _1$ be the disjoint union of $\Gamma '$ and the
 graph on a single vertex $w$ with framing $(-2)$.  The set $S_1$ is
 chosen as $S\cup \{ w\}$.  Simple computation shows that $\mubar
 (Y_{\Gamma _1}, \s _1)=\mubar (Y_{\Gamma }, \s )+1$.  In the surgery
 picture for $Y_{\Gamma _1}$ resulting from the plumbing let $K$
 denote the unknot linking the unknot corresponding to $v\in \Gamma$
 chosen above and the new $(-2)$--framed circle (corresponding to $w$)
 once.
\end{itemize}

Attach a 4--dimensional 2--handle to the 3--manifold $Y_{\Gamma _1}$
along $K$ with framing $(-1)$. The resulting coboridsm will be denoted
by $X$.
\begin{lem}\label{l:topi}
The result of the above surgery is $Y_{\Gamma}$, and the spin
structure $\s _1$ on $Y_{\Gamma _1}$ defined by $S_1$ extends as a
spin structure to provide a spin cobordism $(X, \t )$ from $(Y_{\Gamma
_1}, \s _1)$ to $(Y_{\Gamma}, \s)$.
\end{lem}
\begin{proof}
  By sliding $K$ and the handle corresponding to $w$ down, the first
  statement is obvious. The extension follows from the fact that for
  the graph containing $\Gamma _1$ together with $K$, the vertex
  corresponding to $K$ is not in $S_1$.
\end{proof}

Notice that by induction on $m_{\Gamma ,S}$ the statement of the
theorem holds for $\Gamma _1$ and $\Gamma '$, hence we have that
$-4d(Y_{\Gamma _1},\s _1)=\mubar (Y_{\Gamma _1}, \s _1)= \mubar
(Y_{\Gamma }, \s )+1$ and $-4d(Y_{\Gamma '},\s ')=\mubar (Y_{\Gamma '},
\s ')= \mubar (Y_{\Gamma }, \s)$.

If the spin cobordism $(X, \t)$ of Lemma~\ref{l:topi} between
$(Y_{\Gamma _1}, \s _1)$ and $(Y_{\Gamma}, \s )$ induces a nontrivial
map on the Ozsv\'ath--Szab\'o homology groups, we can easily conclude
the argument: since a negative definite spin cobordism with $\chi =1$
and $\sigma =-1$ shifts degree for Ozsv\'ath--Szab\'o homologies by
$\frac{1}{4}$, the unique nontrivial element of $\hf (Y_{\Gamma _1},
\s _1)$ maps to the unique nontrivial element of $\hf (Y_{\Gamma}, \s
)$ of degree $d(Y_{\Gamma _1}, \s _1)+\frac{1}{4}= d(Y_{\Gamma },
\s)$, reducing the proof to elementary arithmetics. The nontriviality
of the map $F_{X, \t}$ is, however, not so obvious. Let us set up the
exact triangle defined by the surgery triple $(Y_{\Gamma _1},
Y_{\Gamma}, Y_{\Gamma ''})$ along the knot $K\subset Y_{\Gamma _1}$:
\[
\begin{graph}(6,2)
\graphlinecolour{1}\grapharrowtype{2}
\textnode {A}(1,1.5){$\hf (Y_{\Gamma_1})$}
\textnode {B}(5, 1.5){$\hf (Y_{\Gamma})$}
\textnode {C}(3, 0){$\hf (Y_{\Gamma ''})$}
\diredge {A}{B}[\graphlinecolour{0}]
\diredge {B}{C}[\graphlinecolour{0}]
\diredge {C}{A}[\graphlinecolour{0}]
\freetext (3,1.8){$F_{X}$}
\freetext (4.6,0.6){$F_U$}
\freetext (1.4,0.6){$F_V$}
\end{graph}
\]
for the identification of the two manifolds $Y_{\Gamma}, Y_{\Gamma''}$
simple Kirby calculus arguments are needed.  Recall that the map $F_X$
is the sum of all $F_{X, \u}$ for $\u\in Spin^c(X)$.

We claim first that $F_X(\s_1)$ has nonzero $\s$--component.  Since
$U$ is not negative definite, the map $F^{\infty }_U$ vanishes, and
since $Y_{\Gamma}$ is an $L$--space, this implies the same for the
maps $F^+_U$ and $F_U$.  In particular, by exactness we get that $F_V$
is injective and $F_X$ is surjective. Suppose that $F_X(\s _1)$ has
zero $\s$--component. Then $F_X(\s _1) = a +\oa$ for some $a \in \hf
(Y_{\Gamma })$, where $a $ is a formal sum of some spin$^c$ structures
on $Y_{\Gamma}$ and $\oa$ denotes the sum of the conjugate spin$^c$
structures, cf. \cite{OSzII}. By surjectivity now there is $x\in \hf
(Y_{\Gamma _1})$ with $F_X(x )=a $, hence $\s _1+x +\ox$ is in the
kernel of $F_X$, so in the image of $F_V$. If $F_V(y )=\s _1+x +\ox$
then the same holds for $\oy$, hence by the injectivity of $F_V$ the
element $y $ satisfies $\oy=y $. In order $F_V(y)$ to have spin
component, $y$ must have a spin component, hence we have found some
spin and spin$^c$ structures $\z\in Spin (Y_{\Gamma ''})$ and $\t ' \in
Spin ^c (V)$ with $F_{V, \t '}(\z)=\s _1$.  By the uniqueness of
extensions this $\z$ must be $\s ''$, and the spin$^c$ cobordism $(V,
\t ')$ connecting $\z =\s ''$ and $\s _1$ must be spin.  Therefore the
grading shift between the elements $\s ''$ and $\s _1$ is
$\frac{1}{4}$. This implies
\begin{equation}\label{e:de}
d(Y_{\Gamma ''}, \s '')+\frac{1}{4}=d(Y_{\Gamma _1}, \s _1).
\end{equation}
Recall that
\begin{equation}\label{e:mu}
\mubar (Y_{\Gamma ''}, \s '')=\mubar (Y_{\Gamma }, \s )=
\mubar (Y_{\Gamma _1}, \s _1)-1.
\end{equation}
Since by induction for the spin 3--manifolds $(Y_{\Gamma _1}, \s _1)$
and $(Y_{\Gamma ''}, \s '')$ the invariant $\mubar$ actually computes
the correction term, that is, $-4d(Y_{\Gamma _1}, \s _1)=\mubar
(Y_{\Gamma _1}, \s _1)$ and $-4d(Y_{\Gamma '}, \s ')=\mubar (Y_{\Gamma
'}, \s ')$, Equations~\eqref{e:de} and \eqref{e:mu} contradict each
other.  Therefore the element $F_{X}(\s _1)$ has nontrivial
$\s$--component, verifying our claim.

The nontriviality of $F_{X}$ between $\s _1$ and $\s$, however,
implies that there is a connecting spin structure $\t$ with $F_{X, \t
}(\s _1)=\s$, cf. \cite[Lemma~3.3]{OSzII}. Consequently the degree
shift given by $F_{X, \t }$ is $\frac{1}{4}$, hence the inductive step
concludes the proof of Theorem~\ref{t:main}.
\end{proof}

\begin{proof}[Proof of Theorem~\ref{t:obst}]
Combining Propositions~\ref{p:van} and \ref{p:ptlan} with the identification 
of Theorem~\ref{t:main} the proof follows at once.
\end{proof}

{\centerline {{\large {A}}CKNOWLEDGEMENTS} 

We would like to thank Stefan Friedl, Josh Greene and an anonymous
referee for helpful comments and corrections. The author was partially
supported by OTKA 49449, by EU Marie Curie TOK program BudAlgGeo and
by the Clay Mathematics Institute.

\end{document}